\documentclass{amsart}

\newtheorem{conj}{Conjecture}

\def\R{{\bf R}}

\newtheorem{corollary}{Corollary}
\newtheorem{lemma}{Lemma}
\renewcommand{\epsilon}{\varepsilon}
\renewcommand{\phi}{\varphi}
\newtheorem{remark}{Remark}
\newtheorem{theorem}{Theorem}
\newtheorem*{theorem*}{Theorem}
\newtheorem{definition}{Definition}

\newtheorem{prop}{Proposition}

\newcommand{\pa}{\partial}

\begin{document}
\parindent0 in
\parskip 1 em
\title{On the DDVV Conjecture and the Comass in Calibrated Geometry (I)}
\author{Timothy Choi and Zhiqin Lu}
\thanks{The second
author is partially supported by  NSF Career award DMS-0347033 and the
Alfred P. Sloan Research Fellowship.}

\email[Timothy Choi, Department of Mathematics, UC Irvine, Irvine, CA 92697, USA]{tchoi@math.uci.edu}

\email[Zhiqin Lu, Department of Mathematics, UC Irvine, Irvine, CA 92697, USA]{zlu@math.uci.edu}

\date{April 9, 2006}

\maketitle

\tableofcontents

\section{Introduction}
Let $M^n$ be an  $n$ dimensional manifold isometrically immersed into the space form $N^{n+m}(c)$ of constant sectional curvature $c$. Define the normalized scalar curvature $\rho$ and $\rho^\perp$ for the tangent bundle and the normal bundle as follows:
\begin{align}
\begin{split}
&\rho=\frac{2}{n(n-1)}\sum_{1=i<j}^n R(e_i, e_j,e_j,e_i),\\
&\rho^{\perp}=\frac{2}{n(n-1)}\left(\sum_{1=i<j}^n\sum_{1=r<s}^m\langle
R^\perp(e_i,e_j)\xi_r,\xi_s\rangle^2\right)^{\frac 12},
\end{split}
\end{align}
where $\{e_1,\cdots, e_n\}$ (resp. $\{\xi_1,\cdots,\xi_m\}$) is an orthonormal basis of the tangent (resp. normal space) at the point $x\in M$, and $R, R^\perp$ are the curvature tensors for the tangent and normal bundles, respectively.

In the study of submanifold theory, De Smet, Dillen, Verstraelen, and Vrancken~\cite{ddvv} made the following {\sl DDVV Conjecture}: 

\begin{conj}\label{cc1} Let $h$ be the second fundamental form, and 
let $H=\frac 1n \,{\rm trace}\, h$ be the mean curvature tensor.  Then 
\[
\rho+\rho^\perp\leq |H|^2+c.
\]
\end{conj}

 A weaker version of the above conjecture,
 \[
 \rho\leq |H|^2+c,
 \]
 was proved in ~\cite{ch1}. An alternate proof is in~\cite{bogd2}.
 
 In~\cite{ddvv}, the authors proved the following
 \begin{theorem}\label{ddvv}
 If $m=2$, then the conjecture is true.
 \end{theorem}
 
 In this  paper, we prove the conjecture in the case $n= 3$, which is the first non-trivial case.   The proof is quite technical and, like ~\cite{gmm,guw},  some non-trivial linear algebra is involved. We also  point out a relationship between the DDVV conjecture and the comass problem in Calibrated Geometry. We believe that the method used here will be useful in Calibrated geometry.
The more general cases of the conjecture will be treated in our second paper~\cite{lu-30}.

 Let $x\in M$ be a fixed point and let $(h_{ij}^r)$ ($i,j=1,\cdots,n$ and $r=1,\cdots,m$) be the coefficients of the second fundamental form under some orthonormal basis. Then by
 Suceav{\u{a}}~\cite{bogd}, or ~\cite{dfv}, Conjecture~\ref{cc1} can be formulated as an inequality with respect to the coefficients $h_{ij}^r$ as follows:
  \begin{align}\label{1a}
\begin{split}
&\sum_{r=1}^m\sum_{1=i<j}^n(h_{ii}^r-h_{jj}^r)^2+2n\sum_{r=1}^m\sum_{1=i<j}^n(h_{ij}^r)^2\\
&\geq 2n\left(\sum_{1=r<s}^m\sum_{1=i<j}^n\left(\sum_{k=1}^n(h_{ik}^rh_{jk}^s-h_{ik}^sh_{jk}^r)\right)^2\right)^{\frac 12}.
\end{split}
\end{align}

Suppose that $A_1,A_2,\cdots,A_m$ are $n\times n$ symmetric real matrices. Let
\[
||A||^2=\sum_{i,j=1}^n a_{ij}^2,
\]
where $(a_{ij})$ are the entries of $A$, 
and let
\[
[A,B]=AB-BA
\]
be the commutator. Then the equation~\eqref{1a}, in terms of matrices, can be formulated as follows
\begin{conj}\label{conj2}
 For $n, m\geq 2$, we have
\begin{equation}\label{conj}
(\sum_{r=1}^m||A_r||^2)^2\geq 2(\sum_{r<s}||[A_r,A_s]||^2).
\end{equation}
Fixing $n,m$, we call the above inequality Conjecture $P(n,m)$.
\end{conj}

\begin{remark}\label{remk}
For derivation of ~\eqref{1a}, see~\cite[Theorem 2]{dfv}. Note that the prototype of the matrices are the traceless part of the second fundamental forms.
\end{remark}

The main result of this paper is:

\begin{theorem}\label{thm1}
Let $A_i$ $(i=1,\cdots,m)$ be  $3\times 3$ symmetric matrices. Then we have
\begin{equation}\label{main}
(\sum_{i=1}^m||A_i||^2)^2\geq 2\sum_{i<j}(||[A_i,A_j]||^2).
\end{equation}
That is, the Conjecture $P(3,m)$ for $m\geq 2$ is true.
\end{theorem}

In the second part of
of the paper, we discussed the cases when the equality is valid. In particular, we classified all minimal $3$-folds such that the equality of Conjecture~\ref{cc1} is valid at any point. Such kind of $3$-folds are a special kind of {\it austere} submanifolds defined in ~\cite{h-l}. Austere $3$-folds were locally classified in~\cite{bryant-1, df-1}, when the ambient space is the Euclidean space or the unit sphere. We provide a similar classification in the case of  ambient space being hyperbolic, using the  similar method of theirs. Restricting to the situations in this paper, we also give a more straightforward way to classify austere $3$-folds. 

{\bf Acknowledgment.} We thank Suceav{\u{a}} for his bringing this very interesting question, as well as many references, into our attention. We also thank him for  many useful comments, which were absorbed into the current version of this paper. Deep thanks to R. Bryant, who introduced the paper~\cite{df-1} and made many useful comments to me; and  to C. Terng, who explained to me, in great details, of the submanifold geometry in Minkowski spaces; and to M. Dajczer and F. Dillen, who pointed out  to the second author,  that when equality is valid the $3$-fold needs not to be minimal.
Without their helps, the paper won't be in its current form.

\section{Comass in Calibrated Geometry}
Before making further analysis of Conjecture~\ref{conj2}, we introduce the concept of  the comass of a $p$ form in Calibrated Geometry (cf. ~\cite{h-l}). 
 
 Consider Euclidean space $\R^n$ with orthogonal basis $e_1,\cdots,e_n$ and dual basis $dx_i=e_i^*$. Let $I=(i_1,\cdots,i_p)$ denote a multi-index with $i_1<\cdots<i_p$. Let
 \[
 \phi=\sum a_I e^*_{i_1}\wedge\cdots\wedge e_{i_p}^*
 \]
 be a $p$-covector (constant-coefficient $p$-form). 
 The comass $||\phi||^*$ of $\phi$ is given by 
 \[
 ||\phi||^*={\rm max} \{\phi(\xi)\mid\xi \text { is a $p$-plane}\}.
 \]
  For a differential form on a Riemannian manifold $M$, its comass $||\phi||^*$ is given by
 \[
 ||\phi||^*=\sup_x\,\{||\phi_x||^*\mid x\in M\}.
 \]

 In ~\cite{gmm}, Gluck, Mackenzie, and Morgan initiated the study of the comass of the
 first  Pontryagin  form on Grassmann manifolds. Later Gu~\cite{guw} generalized the results, and
 Harvey~\cite{harvey} gave a simplified and unified proof.
 
 Their results are listed as follows:
 
 \begin{theorem*} The comass of the first Pontryagin form $\phi$ on the Grassmann manifold $G(n,m)$ is as follows:
 \begin{enumerate}
 \item $||\phi||^*$ is $\sqrt{3/2}$ for $n=3, m=6$, $4/3$ for $n=3, m\geq 7$, and $3/2$ for $n=4, m\geq 8$~\cite{gmm};
 \item $||\phi||^*$ is $3/2$ for $n\geq 4$ or $m\geq 8$~\cite{guw}.
 \end{enumerate}
 \end{theorem*}

The definition of the comass, in the context of the comass of the first Pontryagin form, can be formulated as following problem:

Let $A,B$ be two $m\times n$ matrices. Define
\[
\{AB\}=AB^T-BA^T.
\]
Let
\begin{align}
\begin{split}
&\qquad \phi(A_1\wedge A_2\wedge A_3\wedge A_4)\\
&=-\frac 12 
{\rm tr}\,(\{A_1A_2\}\{A_3A_4\}+\{A_3A_1\}\{A_2A_4\}+\{A_2A_3\}\{A_1A_4\})
\end{split}
\end{align}
for $m\times n$ matrices $A_1,A_2,A_3$, and $A_4$. The comass of $\phi$  is defined to be the maximum of the right-hand side of the  above under the condition that $A_1,A_2,A_3$, and $A_4$ are orthonormal. 

Conjecture~\ref{conj2} is similar to the above comass problem in that both problems are related to the commutator of matrices. In fact, $P(n,3)$ can be reformulated as follows: let $A,B,C$ be $n\times n$ symmetric matrices such that
\[
||A||^2+||B||^2+||C||^2=1.
\]
Then
\[
||[A,B]||^2+||[B,C]||^2+||[C,A]||^2\leq \frac 12.
\]

The major difference between these  two problems is that they have different invariant groups. The comass problem is invariant under $O(n)\times O(4m)$ (see~\cite{gmm} or ~\cite{guw} for details). On the other  hand the invariant group for the DDVV Conjecture is not known. In the next section, we proved that the invariant group contains the group $O(n)\times O(m)$, which allows us to solve the conjecture for $n=3$.

Although our method is quite different from ~\cite{gmm} or~\cite{guw}, we got hints from the proof of both papers. In particular,  we learnt that the larger the invariant group, the more reductions we can do on the matrices.
Thus the authors strongly believe that the results of both problems should be parallel after more information of the invariant group of the DDVV Conjecture is found. It would also  be very interesting to study the submanifolds ``calibrated'' when the equality of the conjecture is valid at any point.

\section{Invariance}

Let $A_1,\cdots,A_m$ be  $n\times n$ symmetric matrices.
Let $G=O(n)\times O(m)$. Then $G$ acts on matrices $(A_1,\cdots,A_m)$ in the following natural way: let $(p,q)\in G$, where $p,q$ are $n\times n$ and $m\times m$ orthogonal matrices, respectively. Let $q=\{q_{ij}\}$. Then
\[
(p,0)\cdot (A_1,\cdots,A_m)=(pA_1p^{-1},\cdots, pA_mp^{-1}),
\]
and 
\[
(0,q)\cdot (A_1,\cdots,A_m)=(\sum_{j=1}^m q_{1j}A_j,\cdots,\sum_{j=1}^mq_{mj}A_{j}).
\]

It is easy to verify the following 

\begin{prop}\label{prop1}
Conjecture $P(n,m)$ is $G$ invariant. That is, in order to prove inequality ~\eqref{conj}
for $(A_1,\cdots,A_m)$, we just need to prove the inequality for any $\gamma\cdot (A_1,\cdots,A_m)$ where $\gamma\in G$.
\end{prop}

\qed

As a consequence of  the above proposition, we have the following interesting

\begin{theorem}\label{thm3}
Let $n\geq 2$ be an integer. If $P(n, \frac 12 n(n-1)+1)$ is true, then $P(n,m)$ is true  for any $m$.
\end{theorem}

{\bf Proof.}  Obviously, if $P(n, \frac 12 n(n-1)+1)$ is true, then $P(n,m)$ for any $m\leq\frac 12 n(n-1)$ is true because we can set the excessive $A_i$'s to be zero. To prove the case when $m>\frac 12 n(n-1)+1$, we use the math induction. By Proposition~\ref{prop1},  we can assume that at least one of the matrix, say $A_1$, is diagonalized. 
Because the off-diagonal  parts of the matrices form a $\frac 12n(n-1)$-dimensional vector space, 
 we can find  a unit vector $(\alpha_2,\cdots,\alpha_m)$ such that
\[
\alpha_2A_2+\cdots+\alpha_mA_m
\]
is diagonalized. Extending the vectors $(1,0,\cdots,0)$ and  $(0,\alpha_2,\cdots,\alpha_m)$ to an $m\times m$ orthogonal matrix $q$, we can assume, without loss of generality, that $A_2$ is diagnolized. In particular, $[A_1,A_2]=0$.

Replacing $A_1,A_2$ by $\cos\alpha \,A_1+\sin\alpha \,A_2, -\sin\alpha\, A_1+\cos\alpha \,A_2$ for suitable $\alpha$, respectively, we can assume that $A_1\perp A_2$. 

Assuming that  $P(n,m-1)$ is true, we have
\begin{equation}\label{am}
(||A_1+A_2||^2+\sum_{i=3}^m||A_i||^2)^2\geq
2\sum_{i=3}^m||[A_1+A_2, A_i]||^2+2\sum_{3\leq i<j}||[A_i,A_j]||^2,
\end{equation}
and
\begin{equation}\label{am-2}
(||A_1-A_2||^2+\sum_{i=3}^m||A_i||^2)^2\geq
2\sum_{i=3}^m||[A_1-A_2, A_i]||^2+2\sum_{3\leq i<j}||[A_i,A_j]||^2.
\end{equation}

Since $||A_1+A_2||^2=||A_1-A_2||^2=||A_1||^2+||A_2||^2$, adding the above two equations, we have
\[
(\sum_{i=1}^m||A_i||^2)^2\geq 2\left(
\sum_{i=3}^m(||[A_1,A_i]||^2+||[A_2,A_i]||^2)+\sum_{3\leq i<j}||[A_i,A_j]||^2\right).
\]

Since $[A_1,A_2]=0$, the theorem is proved.

\qed

\begin{corollary}\label{cor}
If $P(3,4)$ is true, then $P(3,m)$ is true for $m\geq 2$. Moreover, Theorem~\ref{thm1} follows from $P(3,4)$.
\end{corollary}

\qed

\section{Proof of $P(3,3)$}
We first give a proof of the result of De Smet, Dillen, Verstraelen, and Vrancken~\cite{ddvv}. By our setting in the previous section, the result can be formulated as

\begin{prop}\label{prop2}
Let $A,B$ be $n\times n$ symmetric matrices. Then we have
\[
(||A||^2+||B||^2)^2\geq 2||[A,B]||^2.
\]
\end{prop}

{\bf Proof.} By Proposition~\ref{prop1}, we may assume that one of the matrices, say $A$, is diagnolized.  Let
\[
A=
\begin{pmatrix}
\lambda_1\\
&\ddots\\
&&\lambda_n
\end{pmatrix},\qquad
B=\begin{pmatrix}
b_{11}&\cdots&b_{1n}\\
\vdots&&\vdots\\
b_{n1}&\cdots & b_{nn}
\end{pmatrix}
\]
with $b_{ij}=b_{ji}$. Then we have
\[
||[A,B]||^2=2\sum_{i<j}(\lambda_i-\lambda_j)^2 b_{ij}^2.
\]
Obviously, we have
\begin{equation}\label{k1}
(\lambda_i-\lambda_j)^2\leq2\sum_{k=1}^n\lambda_k^2=2||A||^2.
\end{equation}
Thus we have
\[
||[A,B]||^2\leq 4||A||^2\sum_{i<j}b_{ij}^2\leq 2||A||^2\cdot ||B||^2,
\]
and the proposition follows from the above inequality.

\qed

Before proving $P(3,3)$, we establish the following

\begin{lemma}\label{lem1}
Let $p_1,p_2,p_3\geq 0$ be three real numbers. Let $\lambda_1,\lambda_2,\lambda_3$ be  three
other real  numbers such that $\lambda_1+\lambda_2+\lambda_3=0$ and $\lambda^2_1+\lambda_2^2+\lambda_3^2=1$. Then we have
\begin{equation}\label{1}
\sum_{k=1}^3\lambda_k^2p_k\geq\frac 13 \left(p_1+p_2+p_3-\sqrt{(p_1+p_2+p_3)^2-3(p_1p_2+p_2p_3+p_3p_1)}\right),
\end{equation}
and 
\begin{equation}\label{2}
\sum_{k=1}^3\lambda_k^2p_k\leq\frac 23(p_1+p_2+p_3).
\end{equation}
Furthermore, we have
\begin{align}\label{2-5}
\begin{split}
&\sum_{i<j,k\neq i,j}(\lambda_i-\lambda_j)^2p_k\\
\leq 
&p_1+p_2+p_3+\sqrt{(p_1+p_2+p_3)^2-3(p_1p_2+p_2p_3+p_3p_1)}.
\end{split}
\end{align}
\end{lemma}

{\bf Proof.} Under the restrictions on $\lambda_1,\lambda_2,\lambda_3$,  the maximum and minimum  of the function 
\[
f(\lambda_1,\lambda_2,\lambda_3)=\sum_{k=1}^3\lambda_k^2p_k
\]
exist. Let's use the Lagrange multiplier method to find the maximum and minimum values. Consider the function
\[
F=\sum_{k=1}^3\lambda_k^2p_k-\mu_1(\sum_{k=1}^3\lambda_k)-\mu_2
(\sum_{k=1}^3\lambda_k^2-1),
\]
where $\mu_1,\mu_2$ are multipliers.
At the critical points, we have
\begin{equation}\label{3}
2\lambda_kp_k-\mu_1-2\mu_2\lambda_k=0
\end{equation}
for $k=1,2,3$.
Multiplying by $\lambda_k$ in~\eqref{3} and summing up, we have
\[
\sum_{k=1}^3\lambda_k^2p_k=\mu_2,
\]
at the critical points. Thus the maximum and minimum values of $f$ are the values of  $\mu_2$.

We claim that $\mu_2$ satisfies the equation
\begin{equation}\label{4}
(\mu_2-p_1)(\mu_2-p_2)+(\mu_2-p_2)(\mu_2-p_3)+(\mu_2-p_3)(\mu_2-p_1)=0.
\end{equation}
To see this, we first assume that $\mu_2$ is one of $p_k$ $(k=1,2,3)$. Then by~\eqref{3}, $\mu_1=0$. Since at least two of the three $\lambda_k$'s are not zero, there is $l\neq k$ such that $\mu_2=p_l$. Thus~\eqref{4} is satisfied. On the other hand, if $\mu_2\neq p_1,p_2,p_3$, then from~\eqref{3}, we have
\[
2\lambda_k=\frac{\mu_1}{\mu_2-p_k}
\]
for $k=1,2,3$. Thus $\mu_2$ satisfies ~\eqref{4}.

By solving  $\mu_2$ we get the maximum and minimum values of the function $f$, which proves ~\eqref{1} and ~\eqref{2}. Since
\[
(\lambda_i-\lambda_j)^2=2-3\lambda_k^2,
\]
we get~\eqref{2-5}, and the lemma is proved.

\qed

{\bf Proof of $P(3,3)$.} We let the entries of $A,B,C$ be $a_{ij}, b_{ij}$ and $c_{ij}$, respectively. Using Proposition~\ref{prop1},  we can assume that $A$ is diagonalized.  That is, $a_{ij}=0$ if $i\neq j$.

We let $t^2=a_{11}^2+a_{22}^2+a_{33}^2$. Then  since $m=3$,  inequality ~\eqref{main} can be written as

\begin{equation}\label{8}
(t^2+||B||^2+||C||^2)^2\geq 4t^2\underset{\eta_1^2+\eta_2^2+\eta_3^2=1}{\sum_{i<j,}}(\eta_i-\eta_j)^2(b_{ij}^2+c_{ij}^2)+2||[B,C]||^2.
\end{equation}
Let 
\begin{equation}\label{2-6}
f=\eta_1^2(b_{23}^2+c_{23}^2)+\eta_2^2(b_{13}^2+c_{13}^2)+\eta_3^2(b_{12}^2+c_{12}^2).
\end{equation}
Then ~\eqref{8} can be written as 
\begin{equation}\label{33}
t^4+t^2(-2(||B||^2+||C||^2)+4||\mu||^2+12 f)
+(||B||^2+||C||^2)^2-2||[B,C]||^2\geq 0,
\end{equation}
where $\mu=(b_{11}, b_{22}, b_{33}, c_{11}, c_{22}, c_{33})^T$. 

Let 
\begin{align}
\begin{split}
b_{12}=r_{3}\cos\alpha_3, \quad c_{12}=r_{3}\sin\alpha_3,\\
b_{13}=r_{2}\cos\alpha_2, \quad c_{13}=r_{2}\sin\alpha_2,\\
b_{23}=r_{1}\cos\alpha_1, \quad c_{23}=r_{1}\sin\alpha_1,
\end{split}
\end{align}

Using  ~\eqref{1} of Lemma~\ref{lem1}, the above inequality is equivalent to 
\begin{equation}\label{10}
t^4+t^2\left(2||\mu||^2-4\sqrt{m_0}\right)
+(||B||^2+||C||^2)^2-2||[B,C]||^2\geq 0
\end{equation}
for any $t$, where
\[
m_0=(r_1^2+r_2^2+r_3^2)^2-3(r_1^2r_2^2+r_2^2r_3^2+r_3^2r_1^2).
\]

If $2||\mu||^2-4\sqrt{m_0}\geq 0$, then inequality ~\eqref{10} follows from Proposition~\ref{prop2}. If $2||\mu||^2-4\sqrt{m_0}< 0$, by choosing suitable $t$ minimizing the left hand side of ~\eqref{10}, we get 
\begin{equation}\label{11}
(||B||^2+||C||^2)^2-2||[B,C]||^2\geq (2\sqrt{m_0}-||\mu||^2)^2,
\end{equation}
whenever $2\sqrt{m_0}-||\mu||^2\geq 0$.

Let $\xi\in\R^3$ and let $P: \R^6\rightarrow \R^3$ be a linear map. 
We define
\[
\xi=
\begin{pmatrix}
b_{13}c_{23}-c_{13}b_{23}\\
b_{12}c_{23}-c_{12}b_{23}\\
b_{12}c_{13}-c_{12}b_{13}
\end{pmatrix}.
\]
and
\[
P=
\begin{pmatrix}
c_{12}& -c_{12}&0&-b_{12}&b_{12}&0\\
c_{13}&0&-c_{13}&-b_{13}&0&b_{13}\\
0&c_{23}&-c_{23}&0&-b_{23}&b_{23}
\end{pmatrix}
.
\]
Then  we have
\[
||[B,C]||^2=2||P\mu+\xi||^2.
\]

Using the above  terminology and expanding ~\eqref{11}, we get
\begin{equation}\label{12}
||\mu||^2(r_1^2+r_2^2+r_3^2+\sqrt{m_0}-||Px||^2)-2||\mu||\xi^TPx
+3\sigma_0-||\xi||^2\geq 0,
\end{equation}
where $x=\mu/||\mu||$ and $\sigma_0=r_1^2r_2^2+r_2^2r_3^2+r_3^2r_1^2$.

Since~\eqref{11} and~\eqref{12} are equivalent, ~\eqref{12} is valid if $||\mu||^2=2\sqrt{m_0}$.

If 
\[
r_1^2+r_2^2+r_3^2+\sqrt{m_0}-||Px||^2\leq 0,
\]
or if
\[
r_1^2+r_2^2+r_3^2+\sqrt{m_0}-||Px||^2>0,
\]
but
\[
\xi^TPx\geq \sqrt{2}\cdot\sqrt[4]{m_0}(r_1^2+r_2^2+r_3^2+\sqrt{m_0}-||Px||^2),\quad
\text{or}\quad \xi^TPx\leq 0,
\]
then the minimum of the left hand side of ~\eqref{12} is achieved at either $\mu=0$, or $||\mu||^2=2\sqrt{m_0}$. In view of ~\eqref{11} and the equation ~\eqref{a1} below, in either case, the left hand side of ~\eqref{12} is nonnegative. Thus~\eqref{12} is valid in the above two cases. Finally, if 
\[
r_1^2+r_2^2+r_3^2+\sqrt{m_0}-||Px||^2>0,
\]
and
\begin{equation}\label{cond}
0<\xi^TPx<\sqrt{2}\cdot\sqrt[4]{m_0}(r_1^2+r_2^2+r_3^2+\sqrt{m_0}-||Px||^2),
\end{equation}
then~\eqref{12} is valid if
\[
(\xi^TPx)^2\leq (r_1^2+r_2^2+r_3^2+\sqrt{m_0}-||Px||^2)(3\sigma_0-||\xi||^2).
\]

But by~\eqref{cond}, it suffices to prove
\[
(\xi^TPx)\sqrt{2}\cdot\sqrt[4]{m_0}\leq3\sigma_0-||\xi||^2.
\]

Thus $P(3,3)$ follows from the following

\begin{lemma}\label{lem2}
Using the above notations, we have
\begin{equation}\label{a1}
\sigma_0-||\xi||^2\geq 0,
\end{equation}
and 
\begin{equation}\label{a2}
||P^T\xi|| \sqrt{2}\cdot\sqrt[4]{m_0}\leq 2\sigma_0.
\end{equation}
\end{lemma}

{\bf Proof.} Under the polar coordinates, we have
\[
\xi=
\begin{pmatrix}
r_1r_2\sin(\alpha_1-\alpha_2)\\
r_1r_3\sin(\alpha_1-\alpha_3)\\
r_2r_3\sin(\alpha_2-\alpha_3)
\end{pmatrix}.
\]
Thus~\eqref{a1} follows. To prove~\eqref{a2}, we first write

\[
PP^T=
\begin{pmatrix}
2(c_{12}^2+b_{12}^2)&c_{12}c_{13}+b_{12}b_{13}&-c_{12}c_{23}-b_{12}b_{23}\\
c_{12}c_{13}+b_{12}b_{13}&2(c_{13}^2+b_{13}^2)&
c_{13}c_{23}+b_{13}b_{23}\\
-c_{12}c_{23}-b_{12}b_{23}&c_{13}c_{23}+b_{13}b_{23}&2(c_{23}^2+b_{23}^2)
\end{pmatrix}.
\]

In terms of the polar coordinates, we have

\[
PP^T=
\begin{pmatrix}
2r_3^2& r_2r_3\cos(\alpha_2-\alpha_3)&-r_1r_3\cos(\alpha_3-\alpha_1)\\
r_2r_3\cos(\alpha_2-\alpha_3)& 2r_2^2& r_1r_2\cos(\alpha_1-\alpha_2)\\
-r_1r_3\cos(\alpha_3-\alpha_1) &r_1r_2\cos(\alpha_1-\alpha_2)& 2r_1^2
\end{pmatrix},
\]
from which we have 
\begin{equation}\label{bb}
||P^T\xi||^2=3r_1^2r_2^2r_3^2(\sin^2(\alpha_1-\alpha_2)+\sin^2(\alpha_2-\alpha_3)+
\sin^2(\alpha_3-\alpha_1)).
\end{equation}
Using Lagrange multiplier method, we get 
\[
\sin^2(\alpha_1-\alpha_2)+\sin^2(\alpha_2-\alpha_3)+
\sin^2(\alpha_3-\alpha_1)\leq 9/4.
\]
Thus we have
\[
||P^T\xi||\leq\frac{3\sqrt 3}{2}r_1r_2r_3.
\]
From the above inequality and~\eqref{a2}, it suffices to prove that
\begin{equation}\label{a3}
\frac{3\sqrt 3}{2}r_1r_2r_3 \cdot \sqrt{2}\cdot\sqrt[4]{m_0}\leq 2\sigma_0.
\end{equation}
Let $a,b,c$ be positive numbers. Then by expanding the expression, we have
\[
(a+b+c)^4\geq 6(a^2b^2+b^2c^2+c^2a^2)+4(a^3b+ab^3+b^3c+bc^3+c^3a+ca^3).
\]
By using $a^3b+ab^3\geq 2a^2b^2$, etc, we get
\[
(a+b+c)^4\geq 14(a^2b^2+b^2c^2+c^2a^2).
\]

Using the above inequality, we have
\begin{equation}\label{p-q-r}
\sigma_0^4\geq 14r_1^4r_2^4r_3^4(r_1^4+r_2^4+r_3^4)\geq14r_1^4r_2^4r_3^4m_0.
\end{equation}
Thus
\[
2\sigma_0\geq2\sqrt[4]{14}r_1r_2r_3\sqrt[4]{m_0}\geq \frac{3\sqrt 3}{2}r_1r_2r_3 \cdot \sqrt{2}\cdot\sqrt[4]{m_0},
\]
and ~\eqref{a3} is proved.

\qed

\section{Proof of $P(3,4)$.}

Now we begin to prove $P(3,4)$. That is, we want to prove that, for traceless symmetric $3\times 3$ matrices $A,B,C,D$, we have
\begin{align}
\begin{split}
&(||A||^2+||B||^2+||C||^2+||D||^2)^2\\&\geq
2(||[A,B]||^2+||[A,C]||^2+||[A,D]||^2\\
&\quad +||[B,C]||^2+||[B,D]||^2+||[C,D]||^2).
\end{split}
\end{align}

As before, $a_{ij},b_{ij},c_{ij}, d_{ij}$ represent the $(i,j)$-th entries of the matrices $A,B,C,D$, respectively. Before proving   the inequality, we have the following result which gives some reduction of the matrices:

\begin{lemma}\label{lem3}
Without loss of generality, we may assume
\begin{enumerate}
\item $A$ is diagonalized;
\item $d_{13}=d_{12}=0$;
\item $c_{12}=0$.
\end{enumerate}
\end{lemma}

{\bf Proof.} By  Proposition~\ref{prop1}, we may assume that $A$ is diagonalized. Since the off-diagonal parts of a $3\times 3$ matrix form a three dimensional space, we can find real numbers $\alpha_1,\alpha_2$, and $\alpha_3$ with $\alpha^2_1+\alpha^2_2+\alpha^2_3=1$ such that the entries $(1,3)$ and $(1,2)$ of the matrix $\alpha_1B+\alpha_2C+\alpha_3D$ are zero (after a possible permutation). We now extend $(\alpha_1,\alpha_2,\alpha_3)$ to a $3\times 3$ matrix
\[
\begin{pmatrix}
\alpha_1&\alpha_2&\alpha_3\\
\beta_1&\beta_2&\beta_3\\
\gamma_1&\gamma_2&\gamma_3
\end{pmatrix}.
\]

Replacing $B,C,D$ by $\gamma_1B+\gamma_2C+\gamma_3D$, $\beta_1B+\beta_2C+\beta_3D$, $\alpha_1B+\alpha_2C+\alpha_3D$, respectively, we get $d_{13}=d_{12}=0$. Finally, by choosing suitable $\alpha$ and replacing $B,C$ with $\cos\alpha \,B+\sin\alpha \,C$ and $\sin\alpha\, B-\cos\alpha\, C$  respectively, we may assume that $c_{12}=0$.

\qed

The proof of $P(3,4)$ will be similar to that of $P(3,3)$.
As in the proof of $P(3,3)$, if 
we let $t=||A||$ and let $A'=A/t$, then we have
\begin{align}\label{24}
\begin{split}
&(t^2+||B||^2+||C||^2+||D||^2)^2\\
&\geq 2t^2(||[A',B]||^2+||[A',C]||^2+||[A',D]||^2)\\
&\quad 
+2(||[B,C]||^2+||[B,D]||^2+||[C,D]||^2).
\end{split}
\end{align}

According to Lemma~\ref{lem3}, we assume that $A$ is diagonalized, and we have

\[
B=
\begin{pmatrix}
\mu_1&b_3&b_2\\
b_3& \mu_2&b_1\\
b_2&b_1&\mu_3
\end{pmatrix},\quad
C=\begin{pmatrix}
\mu_4&0&c_2\\
0&\mu_5& c_1\\
c_2&c_1&\mu_6
\end{pmatrix},\quad
D=
\begin{pmatrix}
\mu_7&0&0\\0&\mu_8& d_1\\0&d_1&\mu_9
\end{pmatrix}.
\]

Let 
\begin{align}
\begin{split}
&p_1=\sqrt{b_1^2+c_1^2+d_{1}^2},\\
&p_2=\sqrt{b_{2}^2+c_{2}^2},\\
&p_3=|b_{3}|,\\
&\sigma_1=p_1^2p_2^2+p^2_2p^2_3+p^2_3p^2_1,\\
&m_1=(p_1^2+p^2_2+p^2_3)^2-3\sigma_1,\quad {\rm and} \\
&\mu=(\mu_1,\cdots,\mu_9)^T.
\end{split}
\end{align}

By Lemma~\ref{lem1}, we have
\[
||[A',B]||^2+||[A',C]||^2+||[A',D]||^2\leq 2(p^2_1+p^2_2+p^2_3+\sqrt{m_1}).
\]

If $2\sqrt{m_1}-||\mu||^2\leq 0$, then ~\eqref{24} is trivially true. Otherwise, like ~\eqref{11}, ~\eqref{24} is equivalent to the following

\begin{align}\label{29}
\begin{split}
&\sqrt{(||B||^2+||C||^2+||D||^2)^2-2||[B,C]||^2-2||[B,D]||^2-2||[C,D]||^2}\\
&\geq 2\sqrt{m_1}-||\mu||^2.
\end{split}
\end{align}

Assume that
\[
||[B,C]||^2=2||P_3\mu+\xi_3||^2,
||[B,D]||^2=2||P_2\mu+\xi_2||^2,
||[C,D]||^2=2||P_1\mu+\xi_1||^2.
\]
Then we can write out the matrices explicitly as follows:
\[
P_3=\begin{pmatrix}
0&0&0&-b_3&b_3&0&0&0&0\\
c_2&0&-c_2&-b_2&0&b_2&0&0&0\\
0&c_1&-c_1&0&-b_1&b_1&0&0&0
\end{pmatrix},\quad
\xi_3=\begin{pmatrix}
b_2c_1-b_1c_2\\
b_3c_1\\
b_3c_2
\end{pmatrix},
\]

\[
P_2=\begin{pmatrix}
0&0&0&0&0&0&-b_3&b_3&0\\
0&0&0&0&0&0&-b_2&0&b_2\\
0&d_1&-d_1&0&0&0&0&-b_1&b_1
\end{pmatrix},\quad
\xi_2=\begin{pmatrix}
b_2d_1\\
b_3d_1\\
0
\end{pmatrix},
\]

\[
P_1=\begin{pmatrix}
0&0&0&0&0&0&0&0&0\\
0&0&0&0&0&0&-c_2&0&c_2\\
0&0&0&0&d_1&-d_1&0&-c_1&c_1
\end{pmatrix},\quad
\xi_1=\begin{pmatrix}
c_2d_1\\
0\\
0
\end{pmatrix}.
\]

A straightforward computation gives
\begin{align}
\begin{split}
& \xi_3^TP_3=(b_3c_1c_2, b_3 c_1c_2,-2b_3c_1c_2,-2b_2b_3c_1+b_1b_3c_2, \\
&\qquad\quad  b_2b_3c_1-2b_1b_3c_2, b_2b_3c_1+b_1b_3c_2,0,0,0),\\
&\xi_2^TP_2=(0,0,0,0,0,0,-2b_2b_3d_1, b_2b_3d_1, b_2b_3d_1),\\
&\xi_1^TP_1=0,
\end{split}
\end{align}

We have
\[
||\xi_3^TP_3+\xi_2^TP_2+\xi_1^TP_1||^2
=6b_3^2(c_1^2c_2^2+b_2^2d_1^2+b_2^2c_1^2+b_1^2c_2^2-b_1b_2c_2c_2).
\]
Therefore we have
\[
||\xi_3^TP_3+\xi_2^TP_2+\xi_1^TP_1||^2\leq 7(b_1^2+c_1^2+d_1^2)(b_2^2+c_2^2)b_3^2,
\]
from which we have
\begin{equation}\label{pi}
||\xi_3^TP_3+\xi_2^TP_2+\xi_1^TP_1||\leq\sqrt{7}p_1p_2p_3.
\end{equation}

Using the definition of $p_1,p_2,p_3,\sigma_1,m_1$, and $\mu$, ~\eqref{29} is equivalent to 

\[
(2(p_1^2+p_2^2+p_3^2)+||\mu||^2)^2-4\sum_{i=1}^3||P_i\mu+\xi_i||^2
\geq (2\sqrt{m_1}-||\mu^2||)^2 \quad \text{ if } 2\sqrt{m_1}-||\mu^2||>0.
\]

Extending the above expression, we get the following quadratic inequality:

\begin{equation}\label{31-1}
(p_1^2+p_2^2+p_3^2+m_1-\sum_{i=1}^3||P_ix||^2)||\mu||^2-(\sum_{i=1}^3\xi_i^TP_i)x||\mu||+3\sigma_1-\sum_{i=1}^3||\xi_i||^2\geq 0.
\end{equation}

Again, using the same method as in Lemma~\ref{lem2}, we can prove that
\begin{equation}\label{31}
\sum_{i=1}^3||\xi_i||^2\leq\sigma_1.
\end{equation}

From~\eqref{29} and~\eqref{31}, we know that ~\eqref{31-1} is true if  $||\mu||=0$ or $\sqrt 2\cdot\sqrt[4]{m_1}$. If the minimum value of the quadratic is reached by one of the two points, then the inequality is proved. If the minimum of the above is reached by some point between $(0,\sqrt{2}\cdot\sqrt[4]{m_1})$, then we must have
\begin{enumerate}
\item $p_1^2+p_2^2+p_3^2+m_1-\sum_{i=1}^3||P_ix||^2>0$;
\item $0<\sum_{i=1}^3||P_ix||<\sqrt{2}\cdot\sqrt[4]{m_1}(p_1^2+p_2^2+p_3^2+m_1-\sum_{i=1}^3||P_ix||^2)$.
\end{enumerate}

Like in the proof of $P(3,3)$, it suffices to prove that
\[
\sum_{i=1}^3||\xi_i^TP_i||\cdot\sqrt 2\cdot\sqrt[4]{m_1}\leq 2\sigma_1.
\]

By~\eqref{pi}, it suffices to prove that
\begin{equation}\label{op}
\sqrt{7}p_1p_2p_3\cdot\sqrt 2\cdot\sqrt[4]{m_1}\leq 2\sigma_1.
\end{equation}
By~\eqref{p-q-r}, the above inequality is true, and thus $P(3,4)$ is proved.

\qed

{\bf Proof of Theorem~\ref{thm1}.} Since $P(3,4)$ is true, by Corollary~\ref{cor}, $P(3,m)$ is true.

\qed

\section{The equality cases}

We first establish the following
\begin{prop}\label{prop22}
Let $A,B$ be $n\times n$ symmetric matrices. If
\[
(||A||^2+||B||^2)^2=2||[A,B]||^2,
\]
then there is an orthogonal matrix $Q$ and a real number $\lambda$ such that
\[
A=QA'Q^T,\quad B=QB'Q^T,
\]
where 
\[
A'=\begin{pmatrix}\lambda\\&-\lambda\\&&0\\&&&\ddots\\&&&&0\end{pmatrix},\quad
B'=\begin{pmatrix}0&\pm\lambda\\\pm\lambda&0\\&&0\\&&&\ddots\\&&&&0\end{pmatrix}.
\]
\end{prop}

{\bf Proof.} 
Without loss of generality, we assume $A$ is  diagnolized, then we are in the same situation as in Proposition~\ref{prop2}. In order for  the equality of ~\eqref{k1} to be true, we assume that
\[
(\lambda_1-\lambda_2)^2=\sum_{k=1}^n\lambda_k^2.
\]
Let $\lambda=\lambda_1$. Then we must have $\lambda_2=-\lambda$, $\lambda_k=0$ for $k\geq 2$, and  $b_{ij}=0$ unless $\{i,j\}=\{1,2\}$. Finally, a straightforward computation gives $b_{12}=\pm\lambda$.

\qed

The following theorem shows that the equality cases of $P(3,m)$ is quite restrictive.

\begin{theorem}\label{thm4}
Let $A_i$ $(i=1,\cdots,m)$ be $3\times 3$ traceless matraces. If 
\begin{equation}\label{equal}
(\sum_{i=1}^m||A_i||^2)^2=2(\sum_{i<j}||[A_i,A_j]||^2),
\end{equation}
then up to an element $\gamma\in G$, we have
$A_i=0$ for $i>2$ 
\[
A_1=\begin{pmatrix}\lambda\\&-\lambda\\&&0\end{pmatrix},\quad {\rm and}\quad
A_2=\begin{pmatrix}0&\pm\lambda\\\pm\lambda &0\\&&0\end{pmatrix}.
\]
\end{theorem}

{\bf Proof.} By Propositon~\ref{prop1}, up to an element of $G$, we have
$A_i=0$ for $i>4$. When the equality holds, the equality in~\eqref{29} is reached:
\begin{align}
\begin{split}
&\sqrt{(||B||^2+||C||^2+||D||^2)^2-2||[B,C]||^2-2||[B,D]||^2-2||[C,D]||^2}\\
&= 2\sqrt{m_1}-||\mu||^2.
\end{split}
\end{align}

We claim that one of the matrices  $A,B,C,D$ must be zero (up to a $G$ action). In fact, if $2\sqrt{m_1}-||\mu||^2=0$, then $A=0$. So the claim is valid. If $\mu=0$, then from ~\eqref{31-1} and~\eqref{31}, $\sigma_1=0$. So two of $p_1,p_2,p_3$ are zero. Combining with $\mu=0$, we see that (1). if $p_1=0$, then $D=0$; (2). if  $p_2=p_3=0$, then $B, C, D$ are proportional. Therefore after a $G$ action, one of the matrices  $B,C,D$ must be zero. Finally, if $0<||\mu||^2<2\sqrt{m_1}$, then from~\eqref{op}, again $\sigma_1=0$, and we are in the same situation as in the $\mu=0$ case. 

Since one of the four  matrices $A,B,C,D$ must be zero, we get the equality case of ~\eqref{11}:
\begin{equation}
(||B||^2+||C||^2)^2-2||[B,C]||^2= (2\sqrt{m_0}-||\mu||^2)^2.
\end{equation}

We essentially repeat the previous proof by claiming one of the two matrices $B,C$ must be zero. In fact, if $2\sqrt{m_0}-||\mu||^2=0$, then $A=0$. If $0\leq||\mu||^2<2\sqrt{m_0}$, then  $\sigma_0$ defined in ~\eqref{a2} must be zero. Therefore two of the $r_1,r_2,r_3$ must be zero and $B,C$ must be proportional. Thus after a $G$ action, one of the matrices  $B,C$ must be zero.

In summary, from the above argument, we know that after a $G$ action, only two of the $m$ matrices $A_1,\cdots,A_m$ may be non-zero. Using Proposition~\ref{prop22}, we get the conclusion of the theorem.

\qed

Now we characterize submanifolds $M^3$ of space form $N^{3+m}(c)$ such that $\rho+\rho^\perp=|H|^2+c$ at every points.  By Remark~\ref{remk} and Theorem~\ref{thm4}, we have the following

\begin{theorem}\label{thm5}
Let $M^3$ be a submanifold of the space form $N^{3+m}(c)$ such that
\begin{equation}\label{austere}
\rho+\rho^\perp=|H|^2+c.
\end{equation}
Let the  local orthogonal  frames  of $TM$ be $e_0, e_1,e_2$ and the local orthogonal frames of $T^\perp M$ be $\xi_0,\cdots,\xi_{m-1}$. Let    the matrices of the second fundamental form corresponding to $\xi_i$ be  $\tilde A_i$ for $0\leq i\leq m-1$. Then $\tilde A_i=0$ if $i>2$, and 
\[
\tilde A_0=
\begin{pmatrix}\lambda_0\\&\lambda_0\\&&\lambda_0\end{pmatrix},
\quad
\tilde A_1=
\begin{pmatrix}\lambda_1+\mu\\&\lambda_1-\mu\\&&\lambda_1\end{pmatrix},
\quad
\tilde A_2=
\begin{pmatrix}\lambda_2&\pm\mu\\\pm\mu&\lambda_2\\&&\lambda_2\end{pmatrix},
\]
where $\lambda_0,\lambda_1,\lambda_2$ and $\mu$ are local functions on $M$.
\end{theorem}

\qed

One can easily find  a large classes  of $3$-folds satisfying Theorem~\ref{thm5}. One important class of such $3$-folds were given in~\cite{dfv}, including the totally umbilical submanifolds.

On the other hand, even the classification of {\it minimal} $3$-folds satisfying~\eqref{austere}
will be very interesting and fruitful. So for the rest of the paper, we will only discuss the cases when $M$  is  minimal. The classification is related to the so-called  {\it austere} submanifolds. At the end of this section, we will give the definition. If $M$ is an austere submanifold, then $\lambda_i=0,\, 0\leq i\leq 2$ in Theorem~\ref{thm5}.

Let $e_0,\cdots,e_{2+r}$ be  orthonormal frame fields and 
let $\omega_0,\cdots,\omega_{2+r}$ be   orthonormal coframe fields  of the space form
 $N^{3+r}(c)$. The  Cartan structure equations are as follows
\begin{eqnarray}
&&d\omega_I=-\omega_{IK}\wedge \omega_K;\label{lu-01}\\
&&d\omega_{IJ}=-\omega_{IK}\wedge\omega_{KJ}+c\omega_I\wedge\omega_J,\label{lu-02}
\end{eqnarray}
where the upper case Roman letters $I,J,K,\cdots$ range from $0$ through $2+r$, and $\omega_{IJ}=-\omega_{JI}$. 

We only need to consider the cases $c=0,1,-1$ after rescaling. We use the following standard notations: let  $N^{3+r}(0)=\R^{3+r}$, $N^{3+r}(1)=S^{3+r}\subset\R^{4+r}$, and $N^{3+r}(-1)=H^{3+r}\subset\R^{3+r,1}$, where $\R^{3+r,1}$ is the Minkowski space endowed with the following metric
\begin{equation}\label{lu-02-02}
dx_1^2+\cdots +dx_{3+r}^2-dx_{4+r}^2,
\end{equation}
and $H^{3+r}$ is defined as the hypersurface $x_1^2+\cdots +x_{3+r}^2-x_{4+r}^2=-1$ of $\R^{3+r,1}$.

Let $M$ be an embedded austere $3$-fold in $N^{3+r}(c)$ (where $c=0,1,-1$) such that $e_0,e_1,e_2$ are the  tangent vector fields and $e_3,\cdots,e_{2+r}$ are normal fields of $M$. 
In what follows, we will use lower case Roman letters for the ``normal'' index range
$3\leq a,b,c\leq 2+r$, except for $i,j,k$, where they range from $0$ through $2$.

Along the submanifold $M$, the  Cartan structure  equations~\eqref{lu-01},~\eqref{lu-02} can be re-written as 
\begin{eqnarray}
&& d\omega_i=-\omega_{ij}\wedge\omega_j;\label{lu-1}\\
&&0=-\omega_{aj}\wedge\omega_j;\label{lu-2}\\
&& d\omega_{ij}=-\omega_{iK}\wedge\omega_{Kj}+c\omega_i\wedge\omega_j;\label{lu-3}\\
&&d\omega_{aj}=-\omega_{aK}\wedge\omega_{Kj}.\label{lu-4}
\end{eqnarray}

By the minimality and  Theorem~\ref{thm5},  we may assume that
\begin{equation}\label{lu-5}
\omega_{0a}=0,\quad 3\leq a\leq 2+r.
\end{equation}

From~\eqref{lu-2} and~\eqref{lu-5}, we have
\begin{equation}\label{lu-5-1}
\omega_{a1}\wedge\omega_1+\omega_{a2}\wedge\omega_2=0.
\end{equation}
Thus $\omega_{a1}$ and $\omega_{a2}$ are linear combination of $\omega_1$ and $\omega_2$. Let 
\[
\pi_a=\omega_{a1}-\sqrt{-1}\omega_{a2},
\]
and let 
\[
\omega=\omega_1+\sqrt{-1}\omega_2.
\]
Then by the minimality of $M$, we have\footnote{The equations~\eqref{lu-6} and~\eqref{lu-8} are essentially from ~\cite[\S 4]{bryant-1}.}
\begin{equation}\label{lu-6}
\pi_a=z_a\omega,
\end{equation}
for some complex function $z_a$.  From~\eqref{lu-4},  we have
\[
0=-\omega_{a1}\wedge\omega_{10}-\omega_{a2}\wedge\omega_{20}.
\]
If for some $a$, $z_a\neq 0$, then
 $\omega_{10}, \omega_{20}$ are the linear combinations of $\omega_{a1}, \omega_{a2}$ and hence $\omega_1,\omega_2$. In particular, we have
\begin{equation}\label{lu-7}
\omega_{10}(e_0)=\omega_{20}(e_0)=0.
\end{equation}
We let
\[
\pi_0=\omega_{01}-\sqrt{-1}\omega_{02},
\]
and assume that
\[
\pi_0=z_0\omega+\bar h\bar \omega
\]
for complex functions $z_0$ and $h$.
Since
\[
d\pi_a=-\omega_{aK}\wedge\pi_K,
\]
 by~\eqref{lu-5},~\eqref{lu-6}, we have
\[
0=d\pi_a\wedge\omega=z_a d\omega\wedge\omega=z_a\bar\pi_0\wedge\omega_0\wedge\omega.
\]
Thus $z_0=0$ and 
\begin{equation}\label{lu-8}
\pi_0=\bar h\bar \omega.
\end{equation}

When $h\neq 0$, we can define the $*$-operator on $\Sigma$ by defining $*\omega_{10}=\omega_{20}$, $*\omega_{20}=-\omega_{10}$, and $*1=\omega_{10}\wedge\omega_{20}$.

The following definition is slightly more general than~\cite{h-l}:
\begin{definition}
Let $M^n\rightarrow N^{n+m}(c)$ be an immersed submanifold of the space form. Let $II$ be the second fundamental form of the submanifold. Let $\nu$ be any normal vector of the submanifold. $M^n$ is called {\it austere}, if for any $\nu$ and any $0\leq k<n/2$, we have
$\sigma_{2k+1}(\nu\cdot II)=0$, where $\sigma_k(A)$ is the $k$-th elementary polynomial of the matrix $A$.
\end{definition}

In particular, if $n=3$, austerity is equivalent to $\det (\nu\cdot II)=0$ and minimality.

\begin{corollary}
If the equality in Conjecture~\ref{cc1} is valid at any point and $M$ is minimal, then $M$ is an austere $3$-fold.
\end{corollary}

\qed

\section{On classification of Austere $3$-submanifolds}

Austere $3$-submanifolds of the Euclidean space and unit sphere were locally classified in ~\cite{bryant-1} and ~\cite{df-1}, respectively. Their ideas can be used in the classification of austere $3$-submanifolds of hyperbolic space. In this section,  we will give the local classification of  the  submanifolds for which ~\eqref{austere} hold in different space forms.

\begin{prop}\label{prop4}
Using the above notations, we have
\[
\nabla _{e_0}e_0=0.
\]
Thus locally an austere $3$-fold is a fiber bundle over a $2$-manifold and the fibers are geodesic lines of the space forms. More precisely, there is a neighborhood $U$ of $\R^2$ and  smooth functions $v, u: U\rightarrow \R^{3+r+|c|}$ such that
\begin{enumerate}
\item If $c=0$, then $M$ can be represented by $v+tu$;
\item If $c\neq 0$, then $M$ can be represented by $\cos t\, v+\sin t\, u$, $||u||=1$, $||v||^2=c$, and  $v\perp u$ with respect to the Riemannian or the Minkowski metric in~\eqref{lu-02-02}, respectively.
\end{enumerate}
\end{prop}

{\bf Proof.} If $z_a\equiv 0$ for all $a$, then $M$ is totally geodesic and in this case the theorem follows easily. On the other hand, by~\eqref{lu-5} and ~\eqref{lu-7}, we have $\nabla _{e_0}e_0=0$ on open set where at least one of $z_a\neq 0$. However, since $M$ is minimal, all $z_a$ must be real analytic. Thus by continuity, the above is true everywhere. The local representations follows from the embeddings of $S^{3+r}$ or $H^{3+r}$ into $\R^{4+r}$ or $\R^{3+r,1}$ respectively.

\qed

Let $\Sigma$ be the surface in $M$ defined by $\{t=0\}$. We define the complex structure $J$ on $\Sigma$ as follows: let $P: T\Sigma\rightarrow{\rm span}\{e_1,e_2\}$ be the orthogonal projection. If $Pe=e_1$, then we define $Je=P^{-1}e_2$. It is not hard to see that the definition only depends on the space $\{e_0\}^\perp$, not on the particular choices of $e_1,e_2$.

Let $\langle\quad\rangle_J$ be the Riemannian metric of $\Sigma$ such that $\{e, Je\}$ forms an orthonormal frame. Let $\nabla, \nabla^J$ be the Levi-Civita connection with respect to the induced metric and $\langle\quad\rangle_J$, respectively. Let $T=\nabla^J-\nabla$ and let $X=T(e,e)+T(Je, Je)$. Then we have the following:

\begin{theorem}\label{thm6}
The necessary and sufficient conditions for $M$ to be an austere $3$-fold are the following
\begin{enumerate}
\item If $c=0$, then
\[
du\in{\rm span}\{u, dv\};
\]
if $c\neq 0$, then
\[
du\in{\rm span}\{u, v, dv\};
\]
\item If $c=0$, then 
\[
\Delta v+Xv\in N_\Sigma M;
\]
if $c\neq 0$, then
\[
\Delta v+Xv+\lambda v\in N_\Sigma M, 
\]
where $\lambda=||e||^2+||Je||^2$, and $N_\Sigma M$ is the normal bundle of  $\Sigma\subset M$, which is a real line bundle.
\end{enumerate}
\end{theorem}

{\bf Proof.} Let $M$ be an austere $3$-fold. Let $II$ be the second fundamental form of $M$ in the space forms.  Since $II(e_0,e_1)=II(e_0,e_2)=0$, we have (1); on the other hand, by minimality, let $e=e_1+k_1e_0, Je=e_2+k_2e_0$. Then we have
\begin{equation}\label{lu-11}
II(e,e)+II(Je, Je)=II(e_1,e_1)+II(e_2,e_2)=0,
\end{equation}
which implies  (2).

Now we assume that (1), (2) are valid. Let $(x,y)$ be a local coordinate system of $\Sigma$. If $u_x, u_y$ are linear combinations of $u, v_x, v_y$ ($u,v,v_x,v_y$, resp.), then for $t$ small enough, $u_x, u_y$ ($-sint\, v_x+\cos t\,u_x, -\sin t\, v_y+\cos t \,u_y$, resp.) are linear combinations of $u, v_x+tu_x, v_y+tu_y$ ($-\sin t\, v+\cos t\,u, \cos t\, v_x+\sin t\, u_x, \cos t\, v_y+\sin t\, u_y$, resp.). Thus in  a neighborhood of $t$, $II(e_0, Y)=0$ for any $Y\in T\Sigma$. By the analyticity (which follows  from the  minimality),  $II(e_0,Y)=0$ is true whenever it is defined.

Let $x_0\in M$ be a fixed point and let $\sigma(t)$ be the geodesic line passing through $x_0$ such that $\sigma'(0)=e_0$. Let
$e_0(t), e_1(t), e_2(t)$ be the parallel translation of the frames $e_0, e_1,e_2$ along the geodesic $\sigma(t)$. Then by the assumption that $\omega_{a0}=0$, we know  $e_0(t), e_1(t), e_2(t)$ are parallel even in 
$N^{3+r}(c)$. Using this fact, we see that the equation~\eqref{lu-11}, if true at one point $x_0$, must be true on the geodesic $\sigma(t)$. Finally, (2) is equivalent to ~\eqref{lu-11}. The theorem is proved.

\qed

The above result gives a good classification of austere $3$-folds.  It is also possible to write out all the integrability conditions through the setting. However, as showed in ~\cite{bryant-1}, there are different integrablility conditions in different cases. For the sake of simplicity, for the rest of the paper, we will only locally classify all the {\it generic} austere $3$-folds. 
The setting below is slightly different from that in Theorem~\ref{thm6}: even though we still represent the $3$-folds as $v+tu$, or $\cos t\, v+\sin t\, u$, in what follows, the points $\{t=0\}$ may be singular points of the $3$-folds.
The classification was done
 in~\cite{bryant-1} for $c=0$~\footnote{In the case of $c=0$, the $3$-folds were locally completely classified. Refer to the paper for details.} and in~\cite{df-1} for $c=1$. 

The following result is  from~\cite[Theorem 4.1]{bryant-1}, which solves the case $c=0$:

\begin{theorem}
If $c=0$,  and if the function $h$ in~\eqref{lu-8} is not zero on $M$. Then locally $M$ can be represented by $v+tu$ for functions $v, u: M\rightarrow \R^{3+r}$. By replacing $v$ with $v-\xi u$ if necessary ($\xi$ being a smooth function of $\Sigma$), we have
\begin{align}
&\Delta u=-2u;\\
& dv=*d\phi\, u-\phi*du
\end{align}
where $\phi$ satisfies $\Delta\phi=-2\phi$. Here $\Delta$ is the Laplacian
with respect to $\langle\quad\rangle_J$.
\end{theorem}

The main result of this section is the following theorem. Note that
 when $c=1$, it is  known to Dajczer and Florit~\cite[Theorem 14]{df-1}.
\begin{theorem}\label{thm8}
Let $c=\pm1$,  and let the imaginary part of $h$ be nonzero on $M$. Let $\tilde M$ be the cone over $M$ in the Euclidean space $\R^{4+r}$ ($\R^{3+r,1}$, resp.). Then there is a densely open set $\tilde M^*$ of $\tilde M$ such that $\tilde M^*\cap M$ is also densely open in $M$. Furthermore, $\tilde M^*$ can locally be represented by the regular points of the following map
\[
v+s_1\frac{\pa g}{\pa x}+s_2\frac{\pa g}{\pa y},
\]
where $(x,y)$ is the local coordinate system of $\Sigma$ and $v,g: \Sigma\rightarrow\R^{4+r}$ (or $\R^{3+r,1}$, resp.) are smooth functions; $s_1,s_2$ are two real parameters. In addition, $g$ satisfies
\[
\Delta g+\tilde X(g)=0,
\]
where $\tilde X$ is a vector field and  there are smooth functions $\theta, \phi$ of $\Sigma$ such that
\[
dv=\theta dg+\phi *dg,
\]
and the functions $\theta,\phi$ satisfy the following integrability  conditions:
\begin{enumerate}
\item $d\theta=*(d\phi-\phi \tilde X^*)$;
\item $\Delta\phi-\tilde X(\phi)-{\rm div} \tilde X\phi=0$.
\end{enumerate}
Here $\tilde X^*$ is the $1$ form on $\Sigma$ which is dual to the vector field $\tilde X$ with respect to the metric $\langle\quad\rangle_J$.
\end{theorem}

The proof of the theorem will be similar to that in~\cite{df-1} in spirit. However, our proof is slightly simpler and we avoid using terms like {\it cross sections} or {\it polar surfaces} in that paper. Since the case of $c=1$ is  known and since the proof of $c=1,-1$ are similar, we only give the proof of the theorem when $c=-1$. Recall that we embed the hyperbolic space form $H^{3+r}$ as a hypersurface in the Minkowski space $R^{3+r,1}$ with the metric in~\eqref{lu-02-02}.

We first establish

\begin{prop} If the imaginary part of $h$ is nonzero, then there is a smooth function $g: \Sigma\rightarrow\R^{3+r,1}$ such that
\begin{equation}\label{int-1}
dg=\psi_1 v+\psi_2 e_0,
\end{equation}
for some one forms $\psi_1,\psi_2$ on $\Sigma$. Note that $e_0=u$, where $u$ is defined in Theorem~\ref{thm6}.
\end{prop}

{\bf Proof.}
The integrability condition of ~\eqref{int-1} is
\[
d\psi_1\, v+d\psi_2 \,e_0-\psi_1\wedge(\omega_0 e_0+\omega_1 e_1+\omega_2 e_2)-\psi_2 \wedge(\omega_{10}e_1+\omega_{20}e_2+\omega_0 f)=0,
\]
which is equivalent to the following
\begin{equation}\label{wen-01}
\left\{
\begin{array}{l}
d\psi_1-\psi_2\wedge\omega_0=0;\\
d\psi_2-\psi_1\wedge\omega_0=0;\\
\psi_1\wedge\omega_1+\psi_2\wedge\omega_{10}=0;\\
\psi_1\wedge\omega_2+\psi_2\wedge\omega_{20}=0.
\end{array}
\right.
\end{equation}

From the last two equations of the~\eqref{wen-01}, we get
\begin{equation}\label{wen-1}
-\bar h\psi_2+\psi_1=\xi\pi_0
\end{equation}
for some complex function $\xi$ on $\Sigma$.  Then  from~\eqref{wen-1} and the definition of the $*$-operator, we have
\begin{equation}\label{wen-2}
*(A\psi_2+\psi_1)=B\psi_2,
\end{equation}
where $h=-A-\sqrt{-1}B$ and $A,B$ are smooth real functions of $\Sigma$.

Since $B\neq 0$ by the assumption, the first two equations of ~\eqref{wen-01} can be written as
\[
\left\{
\begin{array}{l}
d\psi_2=(-A\psi_2-B*\psi_2)\wedge\omega_0;\\
 d*\psi_2=\frac{1}{B}(-dB\wedge *\psi_2-dA\wedge\psi_2+A(A\psi_2+B*\psi_2)\wedge\omega_0-\psi_2\wedge\omega_0).
\end{array}
\right.
\]
Let $\psi$ be a solution of 
\[
d\psi=(-A\psi_2-B*\psi_2)\wedge\omega_0.
\]
Then there is a smooth function $\mu$ such that 
\[
\psi_2=\psi+d\mu.
\]
Insert the above equation into the equation of $d*\psi_2$, we get
\begin{equation}\label{wen-001}
\Delta\mu+F(\psi_2)=0,
\end{equation}
where $F(\psi_2)$ is a function of $\psi_2$ that contains the lower order terms. Thus
~\eqref{wen-001} and hence ~\eqref{wen-01} are solvable locally.

\qed

We define a complex structure on the space ${span}\,\{v, e_0\}$ by
\begin{equation}\label{wen-3}
J(e_0-Av)=-Bv.
\end{equation}
We rewrite~\eqref{int-1} as
\begin{equation}\label{int-2}
dg=(\psi_1+A\psi_2) v+\psi_2(e_0-Av).
\end{equation}

We have
\begin{lemma}\label{alex-4}
On smooth $1$-forms of $\Sigma$, we have
\[
*=-J.
\]
\end{lemma}

{\bf Proof.} From ~\eqref{lu-8},  have
\begin{align}\label{59}
\begin{split}
&\omega_{10}=-A\omega_1+B\omega_2;\\
&\omega_{20}=-B\omega_1-A\omega_2.
\end{split}
\end{align}

We let $e\in T\Sigma$ such that $e=e_1+k_1e_0$. Then by the definition of $J$ on $\Sigma$, $Je=e_2+k_2e_0$.  By~\eqref{lu-8}, we have
\[
J(\omega_{10})(e)=\omega_{10}(Je)=\omega_{10}(e_2)=B.
\]
On the other hand, 
\[
*(\omega_{10})(e)=\omega_{20}(e)=\omega_{20}(e_1)=-B,
\]
and the lemma is proved.

\qed

Using the complex structure $J$ on ${span}\,\{v, e_0\}$, we have
\begin{lemma}\label{alex-5}
Let $Y$ be a smooth vector field of $\Sigma$. Then
\[
dg(JY)=J(dg(Y)).
\]
\end{lemma}

{\bf Proof.} Using~\eqref{int-2}, ~\eqref{wen-3}, and Lemma~\ref{alex-4},  we have
\[
dg(JY)=(\psi_1+A\psi_2)(JY) v+\psi_2 (JY)(e_0-Av)=-B\psi_2(Y) v+\frac{\psi_1+A\psi_2}{B}(Y) (e_0-Av).
\]
On the other side
\[
J dg(Y) =J((\psi_1+A\psi_2)(Y)v+\psi_2(Y)(e_0-Av))=\frac{\psi_1+A\psi_2}{B}(Y)(e_0-Av)-B\psi_2(Y) v.
\]
Thus $dg(JY)=J dg(Y)$ and the lemma is proved.

\qed

{\bf Proof of Theorem~\ref{thm8}.}
We let $e\in T\Sigma$ such that
\[
dg(e)=-a(e_0-Av)
\]
for some smooth function $a$ on $\Sigma$.
Then by Lemma~\ref{alex-5}, we have
\[
dg(Je)=aBv.
\]

Let $e=e_1+k_1e_0$. Then $Je=e_2+k_2e_0$ as in the proof of Lemma~\ref{alex-4}. We compute
the $span\,\{e_1,e_2\}$ component of 
\[
\nabla_e dg(e)+\nabla_{Je} dg(Je).
\]
A straightforward computation shows that the it is equal to  $a(-\nabla_{e_1}{e_0}+Ae_1+Be_2)=0$ by~\eqref{59}. Thus we conclude that
\begin{equation}
\Delta g+\tilde X(g)=0.
\end{equation}

Since $II(e_0,Y)=II(v,Y)=0$ for any $Y\in T\tilde M^*$, generically we have
\[
dv=dg\circ S
\]
for some endermorphism $S: T\Sigma\rightarrow T\Sigma$. From the integrable conditions, we know that $S$ commutes with the second fundamental form $II$. Thus  $S$ is a linear combination of the identity morphism and the complex structure $J$, and there are smooth functions $\theta,\phi$ of $\Sigma$ such that
\begin{equation}\label{alex-100}
dv=\theta dg+\phi*dg.
\end{equation}
Let $dV=\omega_{10}\wedge\omega_{20}$ be the volume form of  $\Sigma$. Then by a straightforward computation we have
\[
\tilde X^*=-*(\iota(\tilde X) dV).
\]
Differentiating~\eqref{alex-100} ,we get
\[
d\theta=*(d\phi-\phi \tilde X^*).
\]
Differentiating the above equation, we get
\[
\Delta\phi-\tilde X(\phi)-{\rm div}\, \tilde X\phi=0.
\]
Finally, we observe that $\tilde M^*$ is defined to be the points of $\tilde M$ such that $\{g_x,g_y, v_x+s_1g_{xx}+s_2g_{xy}, v_y+s_1g_{xy}+s_2g_{yy}$ spans a $4$-dimensional space. By the genericity, $\tilde M^*$ is densely open in $\tilde M$. On the other hand, since $x\mapsto x/||x||$ is a submersion, $\tilde M^*\cap M$ is also densely open in $M$.
The theorem is proved.

\qed

In order to characterizing the submanifolds for which ~\eqref{austere} holds at every point, we make the following definition:

\begin{definition}\label{def2}
We say a surface $\Sigma\rightarrow \R^{3+r} (\R^{4+r}, \R^{3+r,1}, resp.)$ satisfies property (A), if there is a complex structure $J$ and a function $u: \Sigma\rightarrow \R^{3+r} (\R^{4+r}, \R^{3+r,1}, resp.)$ such that
\begin{enumerate}
\item If $c=0$, then
\[
du\in{\rm span}\{u, dv\};
\]
if $c\neq 0$, then
\[
du\in{\rm span}\{u, v, dv\};
\]
\item If $c=0$, then 
\[
\Delta v+Xv\in{\rm span}\{u,dv\};
\]
if $c\neq 0$, then
\[
\Delta v+Xv+\lambda v\in {\rm span}\{u,dv\}, 
\]
where $\lambda=||e||^2+||Je||^2$;
\item There are orthonormal frames $\xi_1,\cdots,\xi_r$ in $\R^{3+r} (S^{3+r}, H^{3+r}, resp.)$ which are perpendicular to ${\rm span}\{u, dv\}$ such that the second fundamental form of $v$ on the $\xi_1,\xi_2$ directions are
\[
\begin{pmatrix}
\lambda\\
&-\lambda
\end{pmatrix},\quad
\begin{pmatrix}
&\pm\lambda\\
\pm\lambda
\end{pmatrix},
\]
and the second fundamental form on the directions $\xi_i$ (for $i>2$) are zero, where $\lambda$ is an analytic function of $\Sigma$.
\end{enumerate}
\end{definition}

Combining Proposition~\ref{prop4} and Theorem~\ref{thm6}, we have
\begin{theorem}
Let $\Sigma$ be a surface of property (A). Then there is an austere $3$-fold $M$ constructed in Propostion~\ref{prop4} satisfying~\eqref{austere}. Conversely, any minimal  $3$-fold satisfying~\eqref{austere}
can be reconstructed using the functions $v,u$ in Definition~\ref{def2}  through Proposition~\ref{prop4} and Theorem~\ref{thm6}.
\end{theorem}

\qed

\bibliographystyle{abbrv} 
\bibliography{pub,unp,2007}

\end{document}